\documentclass[12pt]{article}
\usepackage{amsmath}
\usepackage{amssymb}
\newtheorem{teor}{Theorem}[section]

\newtheorem{remar}{Remark}[section]

\newtheorem{question}{Question}

\newcommand{\fdim}{\hspace*{\fill}$Qed$}
\newcommand{\dimostr}{{\bf Proof: }}

\newcommand{\real}{\Bbb{R}}
\newcommand{\complex}{\Bbb{C}}

\newcommand{\K}{K\"{a}hler}

\makeatletter
\newcommand{\Ker}{\mathop{\operator@font Ker}\nolimits}
\newcommand{\trace}{\mathop{\operator@font trace}\nolimits}
\newcommand{\coKer}{\mathop{\operator@font Coker}\nolimits}
\makeatother

\title{Infinite geodesic rays in the space of \K\  potentials}
\author{Claudio Arezzo \and Gang Tian}
\date{}
\begin{document}
\maketitle

\section{Introduction}

It has been shown by Mabuchi (\cite{ma1}), Semmes (\cite{s1})
and Donaldson (\cite{do2}) that the space of \K\ metrics
cohomologous to a fixed one has a riemannian structure of infinite
dimensional symmetric space of negative curvature.
Being its dimension not finite, the usual properties about existence,
uniqueness and regularity for geodesics do not hold a priori.
They appear crucial in a number of applications,
which go from hamiltonian dynamics, Monge-Amp\'{e}re equations and
existence and uniqueness of \K\ - Einstein (and, more generally,
extremal) metrics. One of the most difficult problems is about
the existence of solutions for this geodesic equation.
The only examples known up to now are only trivial ones, i.e., 
those coming from one-parameter subgroups
of automorphisms of the underlying manifold, some special examples on $S^2$ (\cite{do2}), 
and on toric varieties
(\cite{gu}). In \cite{c}, X.X.Chen has proved the existence of
geodesics in a weak sense. He used this weaker version of geodesics to
deduce the uniqueness of extremal metrics when the first Chern class is non-positive.
Despite all this, a general geometric construction is desired and one needs
deeper understanding how the geodesics are related to
the geometry of the underlying manifold.

In \cite{t3}, the K-stability was introduced through special degenerations to study
the existence of K\"ahler-Einstein metrics with positive scalar curvature.
This paper was partly motivated by exploring the connection between special degenerations
and geodesics in the space of K\"ahler metrics. They can be both used to define certain stability
of underlying manifolds. This will provide us a further understanding the geometric invariant
theory in terms of symplectic geometry of the space of K\"ahler metrics.
The other motivation of this paper is to provide a new and general geometric construction
of nontrivial solutions for the geodesic equation. 

Our method to relate the equation for the geodesics directly to
special degenerations of complex structures of underlying
K\"ahler manifolds. We will make use of the one-parameter group action of automorhisms
on a special degeneration of
complex manifolds, however, generic manifold in the special degeneration
may have only finitely many automorhisms.
We will show that if the degeneration is non trivial, then our
general existence result produce geodesic of infinite length.

The same method used to produce these infinite geodesic rays gives also an
existence result for the initial value
problem for geodesic with analytical initial datum. The only known 
explicit examples
of geodesics, given by Donaldson on $S^2$ (\cite{do2}), show that analyticity
is not necessary. On the other hand, our result, combined with the
openess result proven in \cite{do3}, gives a great number of 
geodesics on any manifold with smooth
initial datum. The problem of finding a necessary and sufficient 
condition on the
initial data for the solution of the gedesic equation remains open.

To apply the study of geodesics to the problem of extremal
metrics it is crucial to study the behaviour of the Mabuchi (or K-) energy
along geodesics. In particular it is important to understand
whether the K-energy blows up on a finite length geodesic going to the
boundary. This and other geometrical important 
functionals
along geodesics have been studied by Chen (\cite{ch2}) and Calabi-Chen 
(\cite{cch}).
It seems reasonable to expect that infinite legth geodesics correspond to some
degeneration of the complex structure of the underlying manifold.
This would fit well in the picture of the relationship between the existence of
extremal metrics and some stability property of the underlying 
algebraic manifold. Also it seems very important to understand which property of
these geodesic rays is related to
the nontriviality of those generalized Futaki invariant (\cite{dt}). Note that
generalized Futaki invariants were used to define the K-stability in
\cite{t3}. 

We will end the paper with a discussion of a list of problems related to 
the geometry
of geodesics in the space of \K\ metrics that we believe should be of 
great interest
for future study.

\section{Brief summary of known results}

In this section we recall the geometric setting for
our analytic problem, and in
doing so we try to follow as close as possible the exposition of
Donaldson in \cite{do2}. We refer the reader also to the excellent survey paper
\cite{ch3}.

To start with, we need a compact \K\ manifold
$M$ of complex dimension $n$,
with a \K\ form $\omega _0$.
Then, the space of
\K\ potentials  of \K\ forms in the same cohomology class is given by

$${\mathcal H} = \{ \phi \in {\mathcal C}^{\infty}(M)\mid
\omega _{\phi} = \omega _0 +i \partial \bar{\partial} \phi
>0\}\,\, .$$

Being ${\mathcal H}$ an open subset of the space of smooth
functions on $M$, it is clear that its tangent space is
${\mathcal C}^{\infty}(M)$. Let us also denote by
$dVol_{\phi} = \frac{1}{n!}\omega _{\phi}^n$ the volume
form induced by each function. Given two smooth
functions
$\psi, \chi \in T_{\phi}{\mathcal H}$ for some point
$\phi
\in {\mathcal H}$, we can define their scalar product by
$$ (\psi, \chi)_{\phi} = \int _M \psi \chi dVol_{\phi}
\,\, .$$

Of course such a norm gives the possibility to calcultate
the lenght of a path, and therefore we can study which
curves in ${\mathcal H}$ are geodesics. A direct
calculation gives the first appearance of the geodesic
equation in its simplest form; a path $\phi(t)\in
{\mathcal H}$ is a geodesic if and only if
\begin{equation}
\label{eq}
\phi ''(t) - \frac{1}{2}\mid \mid \nabla _t \phi '(t)\mid
\mid_{\phi(t)}^2 = 0,
\end{equation}
where $'$ stands for differentiation in $t$ and
$\nabla _t$ denotes the covariant derivative for
the metric $\omega _{\phi (t)}$.

The notion of geodesics $\phi(t)$ with a tangent vector field $\psi (t)$,
allows us to define a connection on the
tangent bundle to ${\mathcal H}$ by the following

\begin{equation}
\label{eq"}
D_t(\psi(t))= \frac{\partial \psi(t)}{\partial t} - \frac{1}{2}
(\nabla_t \psi (t), \nabla_t \phi '(t))_{\phi}
\end{equation}

which is torsion-free and metric-compatible.

Of course real constants act on ${\mathcal H}$, so if we call
${\mathcal H}_0 = {\mathcal H}/{\real}$, we can say that,
being our manifold always assumed to be compact,
${\mathcal H}_0$ is the space of \K\ metrics on cohomologous to
$\omega _0$.
Moreover the natural splitting ${\mathcal H}= {\mathcal H}_0 \times
\real$ is riemannian for the structure just defined.

It has been a remarkable discovery due to Mabuchi, Semmes
and Donaldson, how richly this geodesic equation is
related to the geometry of $M$.
\begin{enumerate}
\item
{\em{Geodesics and complex Monge-Amp\`{e}re equations:}}
suppose $\phi (t), t\in [0,1]$ is indeed a geodesic in
${\mathcal H}$. It is convenient to think of this path as
a function
$\Phi\colon M \times [0,1]\times S^1 \rightarrow {\real}$,
given by $\Phi (x,t,s) = \phi _t (x)$, i.e. independent
of the $S^1$ coordinate $s$. Let $\Omega _0$ be the
pull-back of $\omega _0$ to $M \times [0,1]\times S^1$
under the projection map, and put $\Omega_{\Phi} =
\Omega_0 + \partial \bar{\partial}\Phi$ on $M \times
[0,1]\times S^1$ (note that we are taking on $[0,1]\times
S^1$ the standard complex structure).
Then the result is that $\phi (t)$ is a geodesic iff
\begin{equation}
\label{ma}
\Omega _{\Phi}^{n+1} = 0
\end{equation}
(see Proposition $3$ in \cite{do2}).
\item
{\em{Geodesics and Wess-Zumino-Witten equations:}}
The result above induces one to study the complex
Monge-Amp\`{e}re which seems to have a simpler form
than our original equation (\ref{eq}). An even more tempting
approach is to take a Riemann surface $R$ with boundary,
and look at some special maps from $R$ to ${\mathcal H}$.

Let us fix a a map $\rho \colon {\mathcal
C}^{\infty}(\partial(V\times R))$.
A map $f\colon R \rightarrow {\mathcal H}$ with $\rho$ as boundary
data, and whose associated map is $F\colon M\times R \rightarrow
{\real}$,
is said to  satisfy the Wess-Zumino-Witten equations iff it is a
critical point of the functional whose first variation is
(it is a direct calculation to show that such a formula defines
a functional)
\begin{equation}
\label{wzw}
\delta I_{\rho}(F) = \frac{1}{(n+1)!}\int _{M\times R} \delta F
\,\, \Omega_F^{n+1}\,\, .
\end{equation}

Therefore $f$ satisfies the Weiss-Zumino-Witten equations if and only if
$\Omega_F^{n+1}=0$. (see Proposition $4$ in \cite{do2}).
\item
{\em{Geodesics, foliations and holomorphic curves in the
diffeomorphism group of $M$:}}
Let us go back to the solution $\Omega _{\Phi}$ of the
equation \ref{ma} given by a geodesic. $\Omega _{\Phi}$
is a two-form on $M \times [0,1]\times S^1$ which is
closed,
positive on the $M$-slices, and of type $(1,1)$ (note that
in this section we could take any Riemann surface
instead of the cylinder $[0,1]\times S^1$).

Being $\Omega _{\Phi}$ a solution of (\ref{ma}) tells us
that its null space at each point in nontrivial,  and all
the above conditions imply that such a space at every point
has to be a complex line transverse to the $M$-slices.

A very significant variant of this
interpretation is the following (see \cite{s1}): consider on the space of
maps $Map(M,M)$ the natural complex structure given by
the complex structure on $M$. Then a foliation on $M \times
[0,1]\times S^1$ with leaves transverse to the Riemann
surface, can be thought as a map (after fixing a base
point on the Riemann surface)
$F\colon [0,1]\times S^1 \rightarrow Diff(M)$ given by
projection along the leaves, and the fact the leaves are
complex tells us that $F$ is in fact a holomorphic map.
Being $\Omega _{\Phi}$ of type $(1,1)$, we know that the
image of $F$ is contained in the subspace of $Diff(M)$, ${\mathcal Y}$,
of diffeomorphisms which preserve the type of $\omega_0$.
\item
{\em{Geodesics, $K$-energy and uniqueness of Extremal
\K\ metrics:}}
One of the main reasons which raises interest in the
problem of
existence of geodesics is its connection with the
problem of existence and uniqueness of constant scalar
curvature \K\ metrics in a fixed cohomology class (and more generally
of extremal metrics), a well known question
asked by Calabi. The literature on this problem is too vast
to give here any account of it. We point out
\cite{ca1}, \cite{ca2}, \cite{c}, \cite{lb} and \cite{t3}, and
references therein. Here let us just mention that
Mabuchi in \cite{ma1} and \cite{ma2} has defined a functional
(now known as Mabuchi or
$K$-Energy) whose critical points are the extremal
\K\ metrics. It has been observed (Corollary 11 in \cite{do2}) that
the existence of a {\em{smooth}} geodesic connecting
two extremal metrics implies the existence of a holomorphic
automorphism of $M$ which pulls back one metric on the other,
therefore giving essential uniqueness for such a metric.
Such a result can be interpreted by saying that {\em{the $K$-energy
is convex
along geodesic}}, a fact highly nontrivial and useful
in the whole theory.
It is important to note that in his argument
regularity of the geodesic is crucial. X.Chen (\cite{c}) made a
very significant step towards the understanding of the regularity
question for geodesic, proving that the solution to
equation (\ref{ma}) is always ${\mathcal C}^{1,1}$. Geometrically
this means that singularity might occur for example if
null-directions for the symplectic form start arising, but not
blow-ups. Even if this is not smoothness (which cannot be hoped),
he has been able to
prove that it is enough to conclude uniqueness for the extremal
metric in the case of non-positive first Chern class of the manifold.
\end{enumerate}

\section{Solutions to the complex Monge-Amp\`{e}re}

In this section we prove the existence of solutions of equation (\ref{ma})
with analytic boundary data in a sense described below. We will show that
this give rise to a Cauchy problem for which the standard techniques about
convergence of formal solutions of analytic non-characteristic 
equations can be applied.
The solutions interesting in the problems about geodesics have some special
symmetry as we will explain in the next section.
This would reduce our problem to a genuine codimension one problem. 
Nevertheless
we can give a general existence problem for which the solutions with $S^1$
symmetry can be easily extracted.

The first result gives a partial answer
to the initial value problem for geodesics. We discuss plausible
generalizations at the end of the paper.

\begin{teor}
\label{init}
Let $\omega _M$ be a real analytic \K\ metric on $M$ and let
$\psi _0 \colon M \rightarrow {\real}$ be a real analytic function.
Then there exists $\epsilon >0$ and a unique analytic family
of functions $\phi _t(x) = \phi (x,t)$, $\mid t \mid <\epsilon$ s.t.
$(\omega _M  + \partial \bar{\partial}\phi)^{n+1} = 0$ on $M\times
\Delta _{\epsilon}$ and $\frac{d}{dt}\phi _t (x) = \psi _0 (x)$ and
$\phi _0(x) = 0$, where $\Delta _{\epsilon} = \{ S\in {\complex}\mid\,\,
\mid S\mid <\epsilon\}$.
\end{teor}
\dimostr
In this situation one can apply directly Cauchy-Kowalevski theorem provided
we show that (\ref{ma}) is an analytic non-characteristic equation.

This is easily done by writing locally (\ref{ma}) as
$$
det \left( \begin{array}{cc}
g_{\alpha \bar{\beta}} + \frac{\partial^2\phi}{\partial z_{\alpha}\partial
\bar{z}_{\beta}} & \frac{\partial^2\phi}{\partial S \partial z_{\alpha}} \\
  \frac{\partial^2\phi}{\partial S \partial\bar{ z}_{\beta}} &
  \frac{\partial^2\phi}{\partial S \partial \bar{S}}
\end{array}\right) =0
$$
The only thing to observe is of course the positivity of the coefficient
of $ \frac{\partial^2\phi}{\partial S \partial \bar{S}}$, and that
if $\phi $ depends only on $t = Re(S)$, $\frac{\partial \phi}{\partial S}=
\frac{\partial \phi}{\partial t}=\frac{\partial \phi}{\partial\bar{S}}$,
$\frac{\partial^2 \phi}{\partial S \partial \bar{S}} =
\frac{\partial^2 \phi}{\partial t^2}$.

\fdim

In fact it is interesting to write down the solution as a power series
in $t$. To this aim let us write $\phi (x,t) = \sum_{k\geq 1} \theta_k(x) t^k$
so that $\theta _1 = \psi _0$. By plugging this into (\ref{ma}) we
get a recursive formula to define $\theta _k$ inductively. To write
simpler formulae let us just restrict ourselves to the case of
$dim_{{\complex}} M =1$ and $z$ be a local coordinate on $M$.
Then $(\omega _M  + \partial \bar{\partial}\phi)^{2} = 0$ becomes

$$\theta _2 = \frac{\partial \theta_1}{\partial z}\frac{\partial \theta _1}
{\partial \bar{z}} \frac{dz \wedge d\bar{z} \wedge dS \wedge d\bar{S}}
{\omega_M \wedge dS \wedge d\bar{S}},$$
and for $k>2$

$$k(k-1)\theta _k =
\sum_{\stackrel{j+m = k}{j,m \geq 1}} \{ m(m-1)\theta _m 
\frac{\partial ^2 \theta _j}
{\partial z \partial \bar{z}} +
mj (\frac {\partial \theta _m}{\partial z}
\frac{\partial \theta _j}{\partial \bar{z}} +
  \frac {\partial \theta _j}{\partial z}
\frac{\partial \theta _m}{\partial \bar{z}})+$$

$$ + j(j-1) \theta _j \frac{\partial ^2 \theta _m}
{\partial z \partial \bar{z}}\}
  \frac{dz \wedge d\bar{z} \wedge dS \wedge d\bar{S}}
{\omega_M \wedge dS \wedge d\bar{S}}.$$

Without relying on general theorems one can also check that this 
formal solution
converges locally.

This kind of solution of (\ref{ma}) give short time existence for geodesics
starting at an anlytic \K\ metric with analytic initial tangent vector.

The main existence theorem we prove is the following. It is a priori 
independent
of our questions about geodesics. We will show later how to impose symmetry
conditions on the solutions so to be able to use them to construct geodesic
rays.

\begin{teor}
\label{1}
Let $D$ be a divisor of a \K\ manifold $V$ of complex
dimension $n$, and let
$\omega _D$ be an analytic \K\ form defined on $reg(D)= D\setminus sing(D)$.
Let further $\tilde{\omega}_V$ be any closed analytic $(1,1)$ form 
which extends
$\omega _D$
in a neighborhood $V$ of $reg(D)$.

Then there exist
a neighborhood $\tilde{V}$ of $reg(D)$ in $V$ and a function $\phi \colon
\tilde{V} \rightarrow {\real}$ s.t.
\begin{enumerate}
\item
\label{kk}
$\omega = \tilde{\omega}_V + \partial \bar{\partial} \phi$ is smooth on
$\tilde{V}\setminus sing(D)$;
\item
\label{rest}
$\omega _{\mid _{reg(D)}} = \omega _D$;
\item
\label{eqq}
$\omega ^n = 0$
\end{enumerate}
\end{teor}
\dimostr
The method we want to use to construct the function $\phi$ is
closely related to the one used in \cite{cy} and \cite{ty1}. The idea
is to write in a sufficiently small neighborhood of $reg(D)$
\begin{equation}
\phi = \sum_{m>0}\,\,\, \sum_{i+j=m+1} S^i\bar{S}^j\theta _{ij} +
\bar{S}^i{S}^j\overline{\theta _{ij}}
\end{equation}
where $S$ is a defining section of the line bundle, $L_D$, associated to the
divisor, so that $S^i\bar{S}^j$ is a section of $L_D^i \otimes
\overline{L_D}^j$, and the $\theta _{ij}$s are smooth sections of
$L_D^{-i} \otimes \overline{L_D}^{-j}$. Clearly, since
$\phi$ has to be real valued, $\theta _{ij} =
\overline{\theta _{ji}}$.

Let us now fix an hermitian metric $\mid \mid \,\,\,\,\,\, \mid
\mid$ on the normal bundle to $D$ in $V$ and denote by
$\tilde{\omega}$ its curvature form (which is in general
different from $\omega_D$) and by $D$ its covariant derivative
(there should be no confusion with the divisor...).

We now seek (some) $\theta_{ij}$s for which the equation is
satisfied. In order to do this the first step is to calculate
\begin{eqnarray}
\label{ddbar}
\partial\bar{\partial}{\phi} & = &
\sum_{m>0}\,\,\, \sum_{i+j=m+1}\{
ij(S^i\bar{S}^j\theta _{ij} +\bar{S}^i{S}^j\overline{\theta _{ij}})
\frac{DS\wedge\overline{DS}}{\mid S \mid ^2} + \nonumber \\ \vspace{2mm}
&   & i S^i \bar{D}\theta_{ij} \wedge\bar{S}^j\frac{DS}{S}  +
jS^i \bar{S}^j D\theta_{ij}\wedge \frac{\overline{DS}}{\bar{S}} +
\nonumber \\ \vspace{2mm}
&   & iS^j \bar{S}^i D\overline{\theta_{ij}}\wedge
\frac{\overline{DS}}{\bar{S}} + jS^j \bar{S}^i
\bar{D}\overline{\theta_{ij}}\wedge
\frac{DS}{S} + S^iD\bar{D}\bar{S}^j\theta_{ij} +
\nonumber \\ \vspace{2mm}
&   & S^jD\bar{D}\bar{S}^i\overline{\theta_{ij}} +
S^i\bar{S}^jD\bar{D}\theta_{ij} +
S^j\bar{S}^iD\bar{D}\overline{\theta_{ij}} \}
\end{eqnarray}

Being $S$ a holomorphic section we have
\begin{equation}
\begin{array}{ccl}
\bar{D}D S^j & =  & j S^j\tilde{\omega} \\
\bar{D}D \theta _{ij} & = & - D\bar{D}\theta_{ij} - (i-j) \theta_{ij}
\tilde{\omega},
\end{array}
\end{equation}

which readily implies that
\begin{equation}
\label{qq}
(\partial\bar{\partial}\phi)^l\wedge
\tilde{\omega}_V ^{n-l} = \tilde{\omega}^l\wedge \tilde{\omega}_V^{n-l} +
{\mathcal O}(2m+2),
\end{equation}
for $l\geq 2$.

Given a positive integer $k$, we now look at the equation $\omega^n
=0$ disregarding all the terms vanishing along $D$ of order
bigger than $k$ (we will say ``modulo $S^k$'' from now on). For
example by equation (\ref{ddbar}) and (\ref{qq}),
\begin{enumerate}
\item
$\omega^n =0$ modulo $S$ reduces to
\begin{equation}
\label{11}
2\theta_{11}DS\wedge \overline{DS} \wedge \tilde{\omega}_V ^{n-1} +
\tilde{\omega}_V ^n + \sum_{l=2}^n\binom{n}{l} \tilde{\omega} ^l
\wedge  \tilde{\omega}_V ^{n-l}=0
\end{equation}
\item
\label{22}
$\omega^n =0$ modulo $S^2$ reduces to
\begin{eqnarray}
0 & = &
\{ 2SDS\wedge \bar{D}\theta_{20} + 2\theta_{11}DS\wedge \overline{DS}
+ 2S D\theta_{11}\wedge \overline{DS}+  \nonumber \\
& &
2\bar{S} \bar{D}\theta_{11}\wedge DS +  2\bar{S}
D\theta_{02}\wedge \overline{DS} +
4(S\overline{\theta_{12}}+\bar{S}\theta_{12})DS\wedge
\overline{DS}\}
\wedge \tilde{\omega}_V ^{n-1} + \nonumber \\
    & &
  \tilde{\omega}_V ^n
+
\sum_{l=2}^n\binom{n}{l} \tilde{\omega} ^l
\wedge  \tilde{\omega}_V ^{n-l}
\end{eqnarray}
\end{enumerate}
and so on.

We can use (\ref{11}) to define $\theta _{11}$, which is then smooth since
$ DS\wedge \overline{DS} \wedge \tilde{\omega}_V ^{n-1}>0$ near $D$,  so that
the equation $\omega^n =0$ holds modulo $S$. Let us then define
$\theta_{20} = \theta_{02} = 0$. We can use equation (\ref{22})
as a definition for $\theta_{21}$ and $\theta_{12}$ in terms of
$\theta_{11}$, its derivatives and other already fixed pieces.
If we define also $\theta_{30}= \theta _{03} =0$ we can write more
easily the equation modulo $S^3$ when we see the terms coming from equation
(\ref{qq}) coming in.

$\omega^n =0$ modulo $S^3$ reduces to
\begin{eqnarray}
0 & = & {\mbox{all above terms}} +
[S^2(6\theta _{31} DS \wedge \bar{DS} + 2 D\theta_{21}\wedge
\bar{DS}) +  \\
& &
\bar{S}S (8\theta _{22} DS \wedge \bar{DS} + 4 D\theta _{12}
\wedge \bar{DS} + 4 \bar{D} \theta_{21}\wedge DS + \nonumber \\
  & &
2\theta _{11}
\tilde{\omega} +2 D\bar{D}\theta_{11})+
\bar{S}^2(6\theta_{13}DS\wedge \bar{DS} + 2 \bar{D}{S}\wedge DS)]
\wedge \tilde{\omega}_V ^{n-1} + \nonumber \\
  & & \tilde{\omega}_V ^n
+ 4S\bar{S} (\theta_{11}^2\tilde{\omega} + \theta_{11}D\bar{D}\theta_{11})
\wedge DS\wedge \bar{DS}\wedge\tilde{\omega}_V ^{n-2} + \nonumber \\
& &
\sum_{l=2}^n\binom{n}{l} \tilde{\omega} ^l
\wedge  \tilde{\omega}_V ^{n-l}
\end{eqnarray}

It is clear that this phenomenon happens
modulo $S^k$ for any $k$, i.e. the equation modulo $S^k$ can be used
to define inductively all $\theta_{ij}$s with $i+j = k+1$ in terms of
the $\theta_{ij}$s with $i+j \leq k$ and their derivatives.

Defining $\theta _{ij}$ in this way we get a formal power series
which solves our complex Monge-Amp\`{e}re. We claim that this is
in fact a convergent series. This can be seen directly by noticing
that a fixed $\theta _{ij}$, with $i+j = k+1$, is defined by at most
(since cancellations might occur)
\begin{enumerate}
\item
$8$ terms coming from choosing the highest term in $\partial
\bar{\partial} \phi$;
\item
for $i' < i$, $j'<j$ each term with coefficient $S^{i'}\bar{S}^{j'}$
in the expression on $\partial \bar{\partial} \phi$ gets multiplied
by an appropriate terms in the expression of
$\tilde{\omega}_V ^{n-1}(S,\bar{S})$. So how many of these terms do we get?
For each $m= i'+j'$ the
total number is $\frac{1}{2}\sum_{m<k}m(m+1) < \frac{1}{2}k^2(k+1)$.
\item
the remaining terms coming form equation (\ref{qq}).
  Clearly the number of such terms
is bounded by a polynomial in $k$ too.
\end{enumerate}

Adding up the two  contributions, we see that the number of terms
in the expression defining $\theta_{ij}$ grows polynomially as
$i+j$ grows. It is now straightforward (though very lengthy)
to adapt the classical proof of the convergence of the formal solution
with any analytic initial data. In our case the initial data
$\phi _{\mid_{D}}$ and its first derivatives in $S$ and $\bar{S}$ vanish.

\fdim

\medskip

As we mentioned in the previous sections, to apply the above result
to the question of the existence of geodesics we need to study
the question of dependence of the solution on the $S^1$ parameter
in the situation of algebraic degeneration. We first need to set up
correctly the question we face: let $V$ be in fact a holomorphic
family over the unit disc $\Delta \subset {\complex}$ and let
$\pi \colon V \rightarrow \Delta$ the corresponding submersion.
We want to use Theorem \ref{1} when $D= \pi^{-1}(0)$. Clearly
the $S^1$-action on $\Delta$ gives an $S^1$-action on the normal
bundle of $D$ in $V$. We therefore want to know whether we can
choose the solution to the equation (\ref{ma}) equivariantly.
The following theorem answer affermatively:

\begin{teor}
\label{equiv}
In the situation just described, we can choose $\phi$ to satisfy
the extra property of being $S^1$-invariant.
\end{teor}

\dimostr
Of course the $S^1$-action $\sigma$ on the normal bundle
preserves its fibers, and we can consider a defining section to be
an eigenvector for $\sigma$ in the sense that $\sigma ^*(S) = \alpha
S$, for some constant $\alpha$.
Now, we can certainly choose the first extension $\tilde{\omega}_V$ of
$\omega_D$ to be $S^1$-invariant. Moreover we can put on $L_D$
an $S^1$-invariant hermitian metric so that both its covariant derivative
and its curvature form have this property.
The whole point of the proof is now to look at the
recursive formulae for the
$\theta_{ij}$s. For example,
\begin{enumerate}
\item
$\omega^n =0$ modulo $S$ reduces to
\begin{equation}
\label{111}
2\theta_{11}DS\wedge \overline{DS} \wedge \tilde{\omega}_V ^{n-1} +
\tilde{\omega}_V ^n + \sum_{l=2}^n \binom{n}{l} \tilde{\omega} ^l
\wedge  \tilde{\omega}_V ^{n-l}=0 .
\end{equation}
The invariance of $\tilde{\omega}_V$ and of $D$ implies that
$\theta_{11}$ is invariant too.
\item
\label{222}
$\omega^n =0$ modulo $S^2$ reduces to
\begin{eqnarray}
0 & = &
\{ 2\theta_{11}DS\wedge \overline{DS}
+
\bar{D}\theta_{11} \wedge \bar{S}DS + \nonumber \\
& & S D \theta_{11} \wedge
\overline{DS} +
S D\theta_{11}\wedge \overline{DS} +
\bar{S} \bar{D}\theta_{11}\wedge DS + \nonumber \\
    & & 4(S\overline{\theta_{12}}+\bar{S}\theta_{12})DS\wedge
\overline{DS}\}
\wedge \tilde{\omega}_V ^{n-1} + \tilde{\omega}_V ^n
+ \nonumber \\
  & & \sum_{l=2}^n \binom{n}{l} \tilde{\omega} ^l
\wedge  \tilde{\omega}_V ^{n-l}
\end{eqnarray}

This defines $\theta_{12}$ by
\begin{eqnarray}
\bar{S}(4\theta _{12} DS \wedge \overline{DS} + 2DS\wedge
\bar{D}\theta_{11})\wedge
\tilde{\omega}_D ^{n-1} +\nonumber  \\
\tilde{\omega}_D ^n +
\sum_{l=2}^n \binom{n}{l} \tilde{\omega} ^l
\wedge  \tilde{\omega}_V ^{n-l}=0
\end{eqnarray}

and $\theta_{21}$ by
\begin{eqnarray}
S(4\theta_{21} DS\wedge \overline{DS} + 2 D\theta_{11} \wedge
\overline{DS})\wedge \tilde{\omega}_D ^{n-1} + \nonumber
\\ \tilde{\omega}_D ^n
+ \sum_{l=2}^n \binom{n}{l} \tilde{\omega} ^l
\wedge  \tilde{\omega}_V ^{n-l}
  =0.
\end{eqnarray}
Composing the above expressions with $\sigma$ we see that
$\theta_{21}$ has the variance of $\bar{S}$ and
$\theta_{12}$ has the variance of $S$; therefore $S^2\bar{S}\theta_{21}
+ S\bar{S}^2\theta_{12}$ (which is the term modulo $S^3$ in the
expansion of $\phi$) is indeed $S^1$-invariant.

Once again we leave to the reader the simple observation that this
phenomenon occurs at every step of the definition of the $\theta_{ij}$s,
and therefore the proof is complete.
\end{enumerate}

\fdim

One can also try to find directly $S^1$-invariant solutions by relying on
the classical Cauchy-Kowalevski theorem. This is  again a non-characteristic
probelm as in Theorem \ref{init}. Unfortunately if one writes
the local form of the Monge-Amp\`{e}re in polar coordinates one gets an anlytic
equation away from the divisor (which has become a real hypersurface
for the problem with symmetry). So one has to show that these fake 
singularities
of the coefficients do not interfere with the smoothness of the solution.

\begin{remar}
We wish to stress the fact that the solution to the complex
Monge-Amp\`{e}re given in the two theorems above is  not
unique (the choice of $\tilde{\omega}_V$ and of the terms of the form
$\theta_{i,0}$ and
$\theta_{0,j}$ leaves for example great freedom...).
In fact we have the freedom left by a complex gauge (which appears in the
choice of the section $S$) and the action of the symplectomorphisms.
The choice of $\theta_{i,0}$ and
$\theta_{0,j}$ to vanish seems to play the analogue role of the choice of
Bochner coordinates (i.e. those coordinates for which higher order mixed
derivatives of the \K\ potential vanish) when studying the problem of
approximating a polarized
\K\ metric by Bergmann metrics (\cite{t2}).
\end{remar}

\section{Infinite rays}

The aim of this section is to prove that under some
geometric assumption on $V$, the geodesic rays constructed in
the previous section have infinite length.

Even though probably not the most general situation in which
these geodesic will have infinite length, we will prove this
directly in the setup more natural for prooceding to find the
link of theory developed in this paper and the one of obstructions
to the existence of \K\ - Einstein metrics as studied in \cite{t3}.

Let us then recall that an almost Fano variety is an irreducible,
normal variety, such that for some $m$, the pluri-anticanonical
bundle $K^{-m}_{Y_{reg}}$ extends to an ample line bundle, $L$,
over $Y$. In this situation a choice of a basis of $H^0(Y,L^k)$
defines an embedding of $Y$ into some ${\complex}P^N$. Indeed, this can be
reversed in the following sense: if $Y$ is an irreducible, normal subvariety
in some projective space, which is the limit of a sequence of smooth compact
\K\ manifolds of positive first Chern class in the same projective space,
then $Y$ is an almost Fano variety.

Therefore almost Fano varieties arise naturally as degenerations
of Fano, but we want to look at a particular class of degenerations,
so called {\it special}. First recall that a vector field $v$ on $Y$ is
called \emph{admissible}, if it generates a family of automorphisms $\tau _v
(t)$ of $Y$ s.t. $  \tau _v (t)^*L=L$. A degeneration $\pi \colon V
\rightarrow \Delta$ is then called
\emph{special} if there exists a holomorphic vector field $v$ on $V$ such that
$\pi _*v = -t\frac{\partial}{\partial t}$, and which therefore generate the
one-parameter subgroup $z \rightarrow e^{-t}z$ on $\Delta$.

Special degenerations have the property that $\pi^{-1}(t)$ is biholomorphic to
$\pi^{-1}(s)$ whenever $s$ and $t$ are different from zero. On the other hand
$\pi^{-1}(0)$ might not be biholomorphic to the general fibre. We then call
$M$ a Fano smooth manifold isomorphic to $M_t = \pi^{-1}(t)$, and $Y =
\pi^{-1}(0)$.

Being $V$ special, there exists an embedding of $V$ into some product
${\complex}P^N \times \Delta$, in such a way that $v$ induces a one-parameter
subgroup $G= \{\sigma (z)\}_{z\in {\complex}^*} \subset SL(N+1,{\complex})$,
which then satisfies $\sigma (z)(M) = M_{e^{-z}}$. We denote by
$\omega _{FS}$ the
standard Fubini-Study metric on ${\complex}P^N$.

Theorem \ref{equiv} applied to this situation implies that there exists
on some neighborhood $V$ of $reg(Y)$ a function
$\phi \colon V \rightarrow {\real}$ s.t. $\tilde{\omega}_V +\partial
\bar{\partial}\phi$ satisfies the Monge-Amp\`{e}re and $\phi =\phi(p,
\mid S\mid ^2)$, where $S$ is a defining section of the normal bundle
to $Y$. By fixing $\epsilon_0$ sufficiently small we can find a biholomorphism
$\Sigma \colon M_{\epsilon_{0}}\times [x_0, \infty)\times S^1\rightarrow
\pi ^{-1}(\overline{\Delta_{\epsilon _0}} \setminus \{0\}) \subset V$.

$\Sigma$ is really nothing but ``gluing'' the $\sigma (z)$ together
(so one can check that $x_0 = -log \epsilon _0$ even this won't be necessary
in the sequel).

Clearly $\Sigma ^* (\tilde{\omega}_V +\partial
\bar{\partial}\phi)$ still satisfies the Monge-Amp\`{e}re on
$M_{\epsilon_{0}} \times [x_0, \infty) \times S^1$. Moreover
$\phi = \phi (p, \mid S\mid ^2)$, and $\tilde{\omega}_V$ is $S^1$-invariant
as well, and $\phi(\Sigma(p,x,y))= e^{-x-iy}$. Therefore
  $\Sigma ^* (\tilde{\omega}_V +\partial
\bar{\partial}\phi) = \Omega _0 +\partial \bar{\partial}\Phi$, where
$\Phi = \Phi(x), x\in [x_0,\infty)$, does not depend on $y$, and
hence it defines a genuine geodesic on $M_{\epsilon_{0}}$.

We claim that if the complex structure
jumps in the degeneration, i.e. if
$Y$ is not biholomorphic to $M$, then this geodesic has infinite length.

We note that since the geodesic associated to a solution of the
complex Monge-Amp\`{e}re equation is automatically parametrized by a 
multiple of arc length,
and being our curve parametrized on an infinite segment,
the only possibility to contradict the infiniteness of our ray is that
the curve is indeed constant.
This can be seen explicitely in the following way:
the length of the curve of potentials
between  $x_0 = -log \epsilon$ and $x_1$ as $x_1$ goes to infinity
is given by

$$ L= \sqrt{\int_{x_0} ^{x_1}\int_M\mid \Phi'_{x}\mid^2
dVol_{\Omega_{}}dx} = $$

$$ \sqrt{(x_1-x_0)\int_{M}(\Phi')^2_0 dVol_{\Omega_{x_{0}}}}\,\,\, .$$

Therefore the only way to have a finite length curve is that
the family of potential does not depend on $s$. We now claim that
this would contradict
the jump in complex structure.

We now
assume that the degeneration is non trivial and that each component of $Y$ has
multiplicity one. Here $Y$ is allowed to be non smooth and reducible, even
though we do not have a general existence theorem for the solution of the
Monge-Amp\'{e}re equation in the non smooth case.

Indeed, $\mid\mid \Phi_x\mid\mid_{C^0}$, diverges as $x$ goes to $+\infty$.
We can prove this by contradiction: if $\mid\mid \Phi_x\mid\mid_{C^0}$ were
uniformly bounded, then $\sigma(x)\colon M \rightarrow M_{e^{-x}}$
would converge to
a holomorphic map $\tau \colon M \rightarrow Y$. Let us look at $\tau ^{-1}
(p)$ for any $p\in Y$. If $\tau ^{-1} (p)$ is of complex dimension greater
than zero, then for any small neighborhood $U$ of $p$,
$\sigma(x)(\tau ^{-1}(p))$ is contained in $U$, for any sufficiently small
$t$. This is clearly contradicts maximum principle for holomorphic maps.

Similarly, since each component of $Y$ has multiplicity one, $\tau ^{-1} (p)$
has at most one component. Therefore $\tau$ is a biholomorphism and
therefore the degenerartion is trivial, contradicting our assumption.

We summarize the above discussion in the following

\begin{teor}
If $\pi \colon V \rightarrow \Delta$ is a non trivial special degeneration,
whose central fibre is smooth, then the geodesic
produced by solving the degenerate complex Monge-Amp\`{e}re in Theorem \ref{1}
has infinite length.
\end{teor}

A number of natural questions arise naturally from our results. We wish
to point out explicitely some of them.

\begin{question}
It should be possible to generalize our results to special degenerations
such that the central fibre may have mild singularities. 
Moreover it would be also interesting to see what are the
milder assumptions on the initial data for the Monge-Amp\'{e}re equation to get
the convergence of the formal solution. In the similar setting of constructing
\K\ metrics on Grauert tubes (\cite{bu}, \cite{gs}),
analyticity is essential and cannot be dropped
(we wish to thank Prof. D. Burns for extremely useful discussions on this
problem). In the situation studied in this paper more flexibility is allowed
and we know of many examples of geodesics which do not satisfy these 
assumptions.
It is likely that the right assumption is a compatibility condition between
the intial data on $\phi$ and on the given \K\ form.
\end{question}

\begin{question}
It is tempting to draw an analogy between the problem of approximating \K\
metrics by projectively induced metrics, and the one of approximating
geodesics segments by paths induced by subgroups of linear groups of a fixed
embedding into a projective space. Tian (\cite{t2}) proved that any polarized
\K\ metric is a limit of Bergmann metrics as the degree of the projective
embedding goes to infinity.  It would be very interesting to know whether the
analogue "dynamical" result holds for geodesic segments.
\end{question}
\begin{question}
If the previous problem can be solved, one may ask if there is
an expansion of the approximation in term of geometrical data of
the \K\  manifold. For a fixed \K\ metric this is indeed possible and gives
rise to many interesting results (\cite{do3}, \cite{lu}, \cite{z}).
\end{question}

\begin{question}
It is interesting to associate to each geodesic of
infinite lenght a generalized Futaki invariant. It was done for special
degenrations by Tian and Ding-Tian (\cite{t3}, \cite{dt}). This will allow us to
construct analytic obstructions to the
existence of extremal and \K - Einstein metrics in terms of geodesics.
\end{question}

\begin{question}
In order to understand the relationship between the geodesics studied here
and the existence problem for extremal metrics,
we observe that on the geodesics constructed with our method the time
derivative of the K-energy has a limit as time goes to
infinity (this follows from the fact that the sequence
of \K\ metrics converges to a metric). We believe this
property should characterize infinite-length geodesics
corresponding to degenerations of the complex structure.
If this turns out to be true, they should therefore
be the relevant ones to study the stability of the underlying algebraic
manifold, which is the crucial issue for the existence problem for
extremal metrics.
\end{question}

\noindent{Claudio Arezzo} \newline
Dipartimento di Matematica \newline
Universit\`{a} di Parma \newline
Via M. D'Azeglio 85 \newline
43100, Parma \newline
Italy

\noindent{e-mail: claudio.arezzo@unipr.it}

\vspace{3mm}

\noindent{Gang Tian} \newline
Department of Mathematics \newline
Massachusetts Institute of Technology \newline
77 Massachusetts Avenue\newline
02139 Cambridge, MA \newline
U.S.A.

\noindent{e-mail: tian@math.mit.edu}


\begin{thebibliography}{99}

\bibitem{bu}
Burns, D. {\em Curvatures of Monge-Ampère foliations and parabolic manifolds.}
  Ann. of Math. (2) 115 (1982), no. 2, 349-373.

\bibitem{ca1} Calabi E.,
{\em
Extremal \K\ metrics},
in Seminar on Differential Geometry, Annals of Math. Studies {\bf
102} (1982), 259-290.

\bibitem{ca2} Calabi E.,
{\em
Extremal \K\ metrics, II},
in Differential Geometry and Complex analysis, Lecture Notes in Math.
(1985), 96-114, Springer.

\bibitem{cch}
Calabi E. and Chen X.
{\em The space of \K\ metrics II}, preprint, math.DG/0108162.

\bibitem{c} Chen X.,
{\em The space of \K\ metrics},
J. Differential Geom. 56 (2000), no. 2, 189-234.

\bibitem{ch2}
Chen, X.{\em  On the lower bound of the Mabuchi energy and its application}.
  Internat. Math. Res. Notices 2000, no. 12, 607-623.

  \bibitem{ch3}
Chen, X.{\em  Recent progress in \K\ Geometry}.
Proc. of the Internat. Congress of Mathematicians 2002, Vol. II, 273-282.

\bibitem{cy} Cheng S.Y. and Yau S.T.,
{\em Inequality between Chern numbers of singular \K\ surfaces
and characterization of orbit space of discrete group of $SU(2,1)$},
Contemp. Math. {\bf 49} (1986), 31-43.

\bibitem{dt} Ding W. and Tian G.,
{\em \K -Einstein metrics and the generalized Futaki invariant},
Invent. Math. {\bf 110} (1992), 315-335.

\bibitem{do1} Donaldson S.K.,
{\em Remarks on gauge theory, complex geometry, and $4$-manifold topology},
The Fields Medal volume, World Scientific, 1997.

\bibitem{do2} Donaldson S.K.,
{\em Symmetric spaces, \K\  geometry and Hamiltonian dynamics},
Northern California Symplectic Geometry Seminar, 13--33, Amer. Math. Soc.
Transl. Ser. 2, 196, Amer. Math. Soc., Providence, RI, 1999.

\bibitem{do3} Donaldson S.K., {\em Scalar curvature and projective 
embeddings I},
preprint.

\bibitem{do4} Donaldson S.K.,
{\em Holomorphic disks and the complex Monge-Ampere equation},
to appear on Journal of Symplectic Geometry.

\bibitem{gu} Guan, D.,
{\em On modified Mabuchi functional and Mabuchi moduli space of 
KŠhler metrics on toric bundles},
Math. Res. Lett. {\bf 6} (1999), no. 5-6, 547-555.

\bibitem{gs}
Guillemin V. and  Stenzel M. {\em  Grauert tubes and the
homogeneous Monge-Ampère equation}. J. Differential Geom.{\bf  34}
(1991), no. 2, 561-570.

\bibitem{lb} LeBrun C.,
{\em Polarized $4$-manifolds, extremal \K\ metrics, and
Seiberg-Witten theory},  Math. Res. Lett. {\bf 5} (1995), 653-662.

\bibitem{lu} Z. Lu,
{\em On the lower terms of the asymptotic expansion of Tian--Yau--Zelditch},
Amer. J. Math. 122 (2000),  235-273.

\bibitem{ma1} Mabuchi T.,
{\em K-energy maps integrating Futaki invariants},
Tohoku Math. J. (1986), 575-593.

\bibitem{ma2} Mabuchi T.,
{\em Some symplectic geometry on compact \K\ manifolds},
Osaka J. Math. {\bf 24} (1987), 227-252.

\bibitem{s1} Semmes S.,
{\em Complex Monge-Amp\`{e}re and symplectic manifolds},
Amer. J. Math. {\bf 114} (1992), 495-550.

\bibitem{t1} Tian G.,
{\em On \K\ - Einstein metrics on certain \K\ manifolds with $C_(M)>0$},
Invent. Math. {\bf 89} (1987), 225-246.

\bibitem{t2} Tian G.,
{\em On a set of polarized \K\ metrics on algebraic manifolds},
J. Diff. Geometry 32 (1990),  99-130.

\bibitem{t3} Tian G.,
{\em On \K\ - Einstein metrics with positive scalar curvature},
Invent. Math. {\bf 130} (1997), 1-37.

\bibitem{ty1} Tian G. and Yau S.T.,
{\em Complete \K\ manifolds with zero Ricci curvature},
Invent. Math. {\bf 106} (1991), 27-60.

\bibitem{z} S. Zeldtich,
{\em Szeg\"{o} Kernel and a Theorem of Tian},.Internat. Math. Res. Notices
{\bf 6} (1998), 317-331.

\end{thebibliography}
\end{document}